\documentclass{amsart}

\usepackage{amsmath,amsthm,amsfonts,amssymb}
\date{}

\newcommand{\eps}{\varepsilon}
\newcommand{\ol}{\overline}
\newcommand{\IQ}{{\mathbb{Q}}}
\newcommand{\IL}{{\mathbb{L}}}
\newcommand{\IR}{{\mathbb{R}}}
\newcommand{\IN}{{\mathbb{N}}}
\newcommand{\id}{{\operatorname{id}}}

\newtheorem{theorem}{Theorem}
\newtheorem{corollary}{Corollary}
\theoremstyle{definition}
\newtheorem{remark}{Remark}
\newtheorem{example}{Example}
\newtheorem{proposition}{Proposition}

\begin{document}

\title{The regularity of quotient  paratopological groups}

\author{Taras Banakh}
\address{Department of Mathematics, Ivan Franko Lviv National University,
Universytetska 1, Lviv, 79000, Ukraine}
\email{tbanakh@yahoo.com}

\author{Alex Ravsky}
\address{Department of Functional Analysis, Pidstryhach Institute for Applied Problems of Mechanics and Mathematics
National Academy of Sciences of Ukraine,
Naukova 2-b, Lviv, 79060, Ukraine}
\email{oravsky@mail.ru}

\keywords{paratopological group, quotient paratopological group,
group reflexion, regularity}
\subjclass{22A15, 54H10, 54H11}
\begin{abstract}
Let $H$ be a closed subgroup of a regular abelian paratopological
group $G$. The group reflexion $G^\flat$ of $G$ is the group $G$
endowed with the strongest group topology, weaker that the
original topology of $G$. We show that the quotient $G/H$ is
Hausdorff (and regular) if $H$ is closed (and locally compact) in
$G^\flat$. On the other hand, we construct an example of a regular
abelian paratopological group $G$ containing a closed discrete
subgroup $H$ such that the quotient $G/H$ is Hausdorff but not
regular.\end{abstract}

\maketitle

In this paper we study the properties of the quotients of
paratopological groups by their normal subgroups.

By a paratopological group $G$ we understand a group $G$ endowed
with a topology $\tau$ making the group operation continuous, see \cite{ST}. If,
in addition, the operation of taking inverse is continuous, then
the paratopological group $(G,\tau)$ is a topological group. A
standard example of a paratopological group failing to be a
topological group is the Sorgefrey line $\IL$, that is the real
line $\IR$ endowed with the Sorgefrey topology (generated by the
base consisting of half-intervals $[a,b)$, $a<b$).

Let $(G,\tau)$ be a paratopological group and $H\subset G$ be a
closed normal subgroup of $G$. Then the quotient group $G/H$
endowed with the quotient topology is a paratopological group, see
\cite{Ra}. Like in the case of topological groups, the quotient
homomorphism $\pi:G\to G/H$ is open. If the subgroup $H\subset G$
is compact, then the quotient $G/H$ is Hausdorff (and regular)
provided so is the group $G$, see \cite{Ra}. The compactness of
$H$ in this result cannot be replaced by the local compactness as
the following simple example shows.

\begin{example}\label{ex1} The subgroup $H=\{(-x,x):x\in\IQ\}$ is closed and discrete
in the square $G=\IL^2$ of the Sorgenfrey line $\IL$. Nonetheless,
the quotient group $G/H$ fails to be Hausdorff: for any irrational
$x$ the coset $(-x,x)+H$ cannot be separated from zero $(0,0)+H$.
\end{example}

A necessary and sufficient condition for the quotient $G/H$ to be
Hausdorff is the closedness of $H$ in the topology of group
reflexion $G^\flat$ of $G$.

By {\em the group reflexion} $G^\flat=(G,\tau^\flat)$ of a
paratopological group $(G,\tau)$ we understand the group $G$
endowed with the strongest topology $\tau^\flat\subset\tau$
turning $G$ into a topological group. This topology admits a
categorial description: $\tau^\flat$ is a unique topology on $G$
such that\begin{itemize}
\item $(G,\tau^\flat)$ is a topological group; \item the identity
homomorphism $\id:(G,\tau)\to(G,\tau^\flat)$ is continuous; \item
for each continuous group homomorphism $h:G\to H$ into a
topological group $H$ the homomorphism $h\circ \id^{-1}:G^\flat\to
H$ is continuous.
\end{itemize}
Observe that the group reflexion of the Sorgenfrey line $\IL$ is
the usual real line $\IR$.

For so-called 2-oscillating paratopological groups $(G,\tau)$ the
topology $\tau^\flat$ admits a very simple description: its base
at the origin $e$ of $G$ consists of the sets $UU^{-1}$, where $U$
runs over open neighborhoods of $e$ in $G$. Following \cite{BR} we
define a paratopological group $G$ to be {\em 2-oscillating} if
for each neighborhood $U\subset G$ of the origin $e$ there is
another neighborhood $V\subset G$ of $e$ such that $V^{-1}V\subset
UU^{-1}$. The class of 2-oscillating paratopological groups is
quite wide: it contains all abelian (more generally all nilpotent)
as well as saturated paratopological groups. Following I.Guran we
call a paratopological group {\em saturated} if for each
neighborhood $U$ of the origin in $G$ its inverse $U^{-1}$ has
non-empty interior in $G$.

Given a subset $A$ of a paratopological group $(G,\tau)$ we can
talk of its properties in the topology $\tau^\flat$. In
particular, we shall say that a subset $A\subset G$ is {\em
$\flat$-closed} in $G$ if it is closed in the topology
$\tau^\flat$. Also with help of the group reflexion many helpful
properties of paratopological groups can be defined.

A paratopological group $G$ is called
\begin{itemize}
\item {\em $\flat$-separated} if the topology $\tau^\flat$ is
Hausdorff;
\item {\em $\flat$-regular} if it has a neighborhood
base at the origin, consisting of $\flat$-closed sets;
\item {\em $\flat$-compact} if $G^\flat$ is compact.
\end{itemize}

It is clear that each $\flat$-separated (and $\flat$-regular)
paratopological group is functionally Hausdorff (and regular).
Conversely, each Hausdorff (resp. regular) 2-oscillating group is
$\flat$-separated (resp. $\flat$-regular), see \cite{BR}. On the
other hand, there are examples of (nonabelian) Hausdorff
paratopological groups $G$ which are not $\flat$-separated, see
\cite{Ra}, \cite{BR}. The simplest example of a $\flat$-compact
non-compact paratopological group is the Sorgefrey circle
$\{z\in\mathbb C:|z|=1\}$ endowed with the topology generated by
the base consisting of ``half-intervals"
$\{e^{i\varphi}:\varphi\in[a,b)\}$, $a<b$.

Now we are able to state our principal positive result.

\begin{theorem}\label{main1} Let $H$ be a normal subgroup of a
$\flat$-separated paratopological group $G$. Then the quotient
paratopological group $G/H$ is
\begin{enumerate}
\item $\flat$-separated if and only if $H$ is closed in
$G^\flat$;
\item $\flat$-regular if $G$ is $\flat$-regular  and
the set $H$ is locally compact in $G^\flat$.
\end{enumerate}
\end{theorem}

\begin{proof} Let $\pi:G\to G/H$ denote the quotient homomorphism.

1. If $H$ is closed in $G^\flat$ then $G^\flat/H$ is Hausdorff as
a quotient of a Hausdorff topological group $G^\flat$. Since the
identity homomorphism $G/H\to G^\flat/H$ is continuous, the
paratopological group $G/H$ is $\flat$-separated.

Now assume conversely that the paratopological group $G/H$ is
$\flat$-separated. Since the quotient map $\pi^\flat:G^\flat\to
(G/H)^\flat$ is continuous its kernel $H$ is closed in $G^\flat$.
\smallskip

2. Assume that $G$ is $\flat$-regular and $H$ is locally compact
in $G^\flat$. It follows that $H$ is closed in $G^\flat$ (this so
because the subgroup $H\subset G^\flat$, being locally compact, is
complete). Then there is a closed neighborhood $W_1\subset
G^\flat$ of the neutral element $e$ such that the intersection
$W_1\cap H$ is compact in $G^\flat$. Take any closed neighborhood
$W_2\subset G^\flat$ of $e$ such that $W_2^{-1}W_2\subset W_1$. We
claim that $W_2\cap gH$ is compact for each $g\in G$. This is
trivial if $W_2\cap gH$ is empty. If not, then $gh=w$ for some
$h\in H$ and $w\in W_2$. Hence $W_2\cap gH\subset W_2\cap
wh^{-1}H=W_2\cap wH=w(w^{-1}W_2\cap H)\subset w(W_2^{-1}W_2\cap
H)\subset w(W_1\cap H)$ and the closed subset $W_2\cap gH$ of $G$
lies in the compact subset $w(W_1\cap H)$ of $G$. Consequently,
$W_2\cap gH$ is compact for any $y\in G$. Let $W_3\subset G^\flat$
be a neighborhood of $e$ such that $W_3^{-1}W_3\subset W_2$.

To prove the $\flat$-regularity of the quotient group $G/H$, given
any neighborhood $U\subset G$ of $e$ it suffices to find a
neighborhood $V\subset U$ of $e$ such that $\pi(V)$ is
$\flat$-closed in $G/H$. By the $\flat$-regularity of $G$, we can
find a $\flat$-closed neighborhood $V\subset U\cap W_3$. We claim
that $\pi(V)$ is $\flat$-closed in $G/H$. Since the identity map
$(G/H)^\flat\to G^\flat/H$ is continuous, it suffices to verify
that $\pi(V)$ is closed in the topological group $G^\flat/H$.

Take any point $gH\notin\pi(V)$ of $G^\flat/H$. It follows from
$gH\cap V=\emptyset$ and the compactness of the set $W_2\cap gH$
that there is an open neighborhood $W_4\subset W_3$ of $e$ in
$G^\flat$ such that $W_4(W_2\cap gH)\cap V=\emptyset$. We claim
that $W_4z\cap V=\emptyset$ for any $z\in gH$. Assuming the
converse, find a point $v\in W_4z\cap V$. It follows that $z\notin
W_2$. On the other hand, $z\in W_4^{-1}v\subset W_4^{-1}V\subset
W_2$. This contradiction shows that $W_4gH\cap V=\emptyset$ and thus
$\pi(W_4g)$ is a neighborhood of $gH$ in $G^\flat/H$, disjoint
with $\pi(V)$.
\end{proof}

\begin{corollary} If $H$ is a $\flat$-compact normal subgroup of a $\flat$-regular
paratopological group $G$, then the quotient paratopological group
$G/H$ is $\flat$-regular.
\end{corollary}

\begin{proof} It follows that the identity inclusion
$H^\flat\to G^\flat$ is continuous and thus $H$ is compact in
$G^\flat$. Applying the preceding theorem, we conclude that the
quotient group $G/H$ is $\flat$-regular.
\end{proof}

\begin{remark} It is interesting to compare the latter corollary with a result
of \cite{Ra} asserting that the quotient $G/H$ of a Hausdorff
(regular) paratopological group $G$ by a compact normal subgroup
$H\subset G$ is Hausdorff (regular).
\end{remark}

Since for a 2-oscillating paratopological group $G$ the Hausdorff
property (the regularity) of $G$ is equivalent to the
$\flat$-separatedness (the $\flat$-regularity),
Theorem~\ref{main1} implies

\begin{corollary}\label{cor1} Let $H$ be a normal subgroup of a
Hausdorff 2-oscillating paratopological group $G$. Then the
quotient paratopological group $G/H$ is
\begin{enumerate}
\item Hausdorff if $H$ is closed in
$G^\flat$;
\item regular if $G$ is regular  and
the set $H$ is locally compact in $G^\flat$.
\end{enumerate}
\end{corollary}

Example~\ref{ex1} supplies us with a locally compact closed
subgroup $H$ of a $\flat$-regular paratopological group $G=\IL^2$
such that the quotient $G/H$ is not Hausdorff. Next, we construct
a $\flat$-regular abelian paratopological group $G$ containing a
locally compact $\flat$-closed subgroup $H$ such that the quotient
is Hausdorff but not regular. This will show that in
Theorem~\ref{main1} and Corollary~\ref{cor1} the local compactness
of $H$ in $G^\flat$ cannot be replaced by the local compactness
plus $\flat$-closedness of $H$ in $G$.

Our construction is based on the notion of a {\em cone topology}
(see the paper~\cite{Ra4} of the second author).
Let $G$ be a topological group and $S\subset G$ be a closed
subsemigroup of $G$, containing the neutral element $e\in G$. The
{\em cone topology} $\tau_S$ on $G$ consists of sets $U\subset G$
such that for each $x\in U$ there is an open neighborhood
$W\subset G$ of $e$ such that $x(W\cap S)\subset U$. It is clear
that the group $G$ endowed with the cone topology $\tau_S$ is a
regular paratopological groups and its neighborhood base at $e$
consists of the sets $W\cap S$, where $W$ is a neighborhood of $e$
in $G$. Moreover, the paratopological group $(G,\tau_S)$ is
saturated if $e$ is a cluster point of the interior of $S$ in $G$.
In the latter case the paratopological group $(G,\tau_S)$ is
2-oscillating and thus $\flat$-regular, see \cite[Theorem 3]{BR}.
\smallskip

In the following example using the cone topology we construct a
saturated regular paratopological group $G$ containing a
$\flat$-closed discrete subgroup $H$ with non-regular quotient
$G/H$.

\begin{example} Consider the group $\IQ^3$ endowed with the usual (Euclidean)
topology. A subsemigroup $S$ of $\IQ^3$ is called a {\em cone} in
$\IQ^3$ if $q\cdot \vec x\in S$ for any non-negative $q\in\IQ$ and
any vector $\vec x\in S$.

 Fix a sequence $(z_n)$ of rational numbers such that
$0<\sqrt{2}-z_n<2^{-n}$ for all $n$ and let $S\subset \IQ^3$ be
the smallest closed cone containing the vectors $(1,0,0)$ and
$(\frac1n,1,z_n)$ for all $n$. Let $\tau_S$ be the cone topology
on the group $\IQ^3$ determined by $S$. Since the origin of
$\IQ^3$ is a cluster point of the interior of $S$, the
paratopological group $G=(\IQ^3,\tau_S)$ is saturated and
$\flat$-regular. Moreover, its group reflexion coincides with
$\IQ^3$.

Now consider the $\flat$-closed subgroup $H=\{(0,0,q):q\in\IQ\}$
of the group $G$. Since $H\cap S=\{(0,0,0)\}$, the subgroup $H$ is
discrete (and thus locally compact) in $G$. On the other hand $H$
fails to be locally compact is $\IQ^3$, the group reflexion of
$G$.

We claim that the quotient group $G/H$ is not regular. Let
$\pi:G\to G/H$ denote the quotient homomorphism. We can identify
$G/H$ with $\IQ^2$ endowed with a suitable topology.

Let us show that $(0,1)\notin\pi(S)$. Assuming the converse we
would find $x\in\IQ$ such that $(0,1,x)\in S$. It follows from the
definition of $S$ that $x\ge0$ and there is a sequence $(\vec
x_i)$ converging to $(0,1,x)$ such that

$$\vec x_i=\sum_n\lambda_{i,n}(n^{-1},1,z_n)+\lambda_i(1,0,0)$$
where all $\lambda_i,\lambda_{in}\ge 0$ and almost all of them
vanish. Taking into account that $\{\vec x_i\}$ converges to
$(0,1,x)$ we conclude that
\begin{itemize}
\item $\lambda_i\to0$ as $i\to\infty$;
\item
$\lambda_{in}\underset{i\to\infty}\longrightarrow0$ for every $n$;
\item $\sum_n\lambda_{in}$ tends to $1$ as $i\to\infty$.
\end{itemize}

Let $\eps>0$. Then

$\exists N_1(\forall n>N_1)\{|z_n-\sqrt{2}|<\eps\}$,

$\exists N_2(\forall i>N_2)(\forall n\le
N_1)\{\lambda_{in}<\eps/N_1\}$ and

$\exists N_3(\forall i>N_3)(|\sum \lambda_{in}-1|<\eps\}$.

Put $N=\max\{N_2,N_3\}$. Let $i>N$. Then

$$|\sqrt{2}-\sum_n\lambda_{in}z_n|\le
|\sqrt{2}-\sum_n\lambda_{in}\sqrt{2}|+ |\sum_{n\le
N_1}\lambda_{in}(\sqrt{2}-z_n)|+|\sum_{n>N_1}\lambda_{in}(\sqrt{2}-z_n)|\le$$
$$\eps\sqrt{2}+\eps+\sum_{n>N_1}\lambda_{in}\eps\le
\eps(\sqrt{2}+1+1+\eps).$$

So $x=\sqrt{2}$ which is impossible. This contradiction shows that
$(0,1)\notin\pi(S)$ and thus $(0,\frac1n)\notin\pi(S)$ for all
$n\in\IN$ (since $S$ is a cone).

It remains to prove that for each neighborhood $V\subset \IQ^3$ of
the origin we get $\overline{\pi(V\cap S)}\not\subset\pi(S)$,
where the closure is taken in $G/H$. This will follow as soon as
we show that $(0,\frac1m)\in\overline{\pi(V\cap S)}$ for some $m$.
Since $V$ is a (usual) neighborhood of $(0,0,0)$ in $\IQ^3$, there
is $m\in\IN$ such that $\frac1m(\frac1n,1,z_n)\in V$ for all
$n\in\IN$. Then $\frac1m(\frac1n,1)\in\pi(V\cap S)$ for all
$n\in\IN$. Observe that the sequence $\{(\frac1{nm},\frac1m)\}_n$
converges to $(0,\frac1m)$ in $G/H$ since for each neighborhood
$W\subset\IQ^3$ of $(0,0,0)$ the difference
$(\frac1{nm},\frac1m)-(0,\frac1m)=(\frac1{nm},0)$ belongs to
$\pi(W\cap S)$ for all sufficiently large $n$. Therefore
$(0,\frac1m)\in\overline{\pi(V\cap S)}\not\subset
\pi(S)\not\ni(0,\frac1m)$, which means that $G/H$ is not regular.
\end{example}

As we understood, in the submitted version of the paper~\cite{XieLiTu}
Li-Hong Xie, Piyu Li, and Jin-Ji Tu proved that if $\mathcal P$ is one of the following properties
$\{T_1,T_2,T_3,$ $regular\}$, a paratopological group $G$ has the property $\mathcal P$,
and $H$ is a compact normal semigroup of the group $G$ then the quotient group $G/H$
has the property $\mathcal P$ too. But the case when $\mathcal P=T_0$ was remarked as unknown.
We fill this gap here.

\begin{proposition} Let $H$ be a compact normal subgrop of a $T_0$ paratopological group $G$.
Then the quotient group $G/H$ is $T_0$ too.
\end{proposition}
\begin{proof} Let $\mathcal B$ be the family of all open neighborhoods of the unit
of the group $G$ and $\mathcal B'$ be the family of all open neighborhoods of the unit
of the group $G/H$.
Let $\pi:G\to G/H$ be the quotient map.
Let $S=\bigcap_{U\in\mathcal B} U$ and $S'=\bigcap_{U'\in\mathcal B'} U'$.
Then $S'\subset \bigcap_{U\in\mathcal B}\pi(UH)\subset \pi(\bigcap_{U\in\mathcal B}UH)$.
Let $x\in \bigcap_{U\in\mathcal B}UH$ be an arbitrary point and $U\in\mathcal B$ be an
arbitrary neighborhood. There exists a neighborhood $V\in\mathcal B$ such that $V^2\subset U$.
Then  $U^{-1}x\supset \ol{V^{-1}x}\supset V^{-1}x$. So $U^{-1}x\cap H\supset
\ol{V^{-1}x}\cap H\ne\emptyset$. Since the set $H$ is compact there exists
point $y\in \bigcap_{U\in\mathcal B}(\ol{U^{-1}x}\cap H)=\bigcap_{U\in\mathcal B}(U^{-1}x\cap H)$.
So $x\in Sy\subset SH$. Hence $S'\subset\pi(SH)$ and
$S'\cap S'^{-1}\subset\pi(SH)\cap \pi(SH)^{-1} \subset\pi(SH\cap S^{-1}H)$.
Let $x\in SH\cap S^{-1}H$ be an arbitrary point. Then there exist elements $s_1,s_2\in S$ and
$h_1,h_2\in H$ such that $x=s_1h_1=s_2^{-1}h_2$. Then $s_2s_1=h_2h_1^{-1}\in S\cap H$.
But since $H$ is a compact paratopological group, by Lemma 5.4 from~\cite{Ra3}, $H$ is
a topological group. Since $H$ is a $T_0$ topological group the space $H$ is $T_1$ (in fact,
$T_{31/2}$), so $H\cap S=H\cap \bigcap_{U\in\mathcal B} U=\{e\}$. Thus $s_2s_1=h_2h_1^{-1}=e$,
so $s_2=s_1^{-1}$ and $h_2=h_1$. Then $xh^{-1}\in S\cap S^{-1}=\{e\}$. Hence $x\in H$.
At last, $S'\cap S'^{-1}\subset \pi(SH\cap S^{-1}H)\subset \pi(H)=\{e\}$ and thus
the group $G/H$ is $T_0$.

\end{proof}


\begin{thebibliography}{}
\bibitem[BR]{BR}
T.~Banakh, A.~Ravsky. Oscillator topologies on a paratopological
group and related number invariants, // {\em Algebraical
Structures and their Applications}, Kyiv: Inst. Mat. NANU, (2002)
140-153.

\bibitem[Ka]{Ka}
M. Kat\u{e}tov, {\em On H-closed extensions of topological spaces},
~\u{C}asopis P\u{e}st. Mat. Fys. {\bf 72} (1947), 17--32.

\bibitem[Ra]{Ra}
A.~Ravsky, Paratopological groups I, {\em Matematychni Studii},
{\bf 16}:1 (2001), 37--48. http://matstud.org.ua/texts/2001/16$\_$1/37$\_$48.pdf

\bibitem[Ra2]{Ra2}
A. Ravsky, {\em Paratopological groups II},
Matematychni Studii {\bf 17}, No. 1 (2002), 93--101. http://matstud.org.ua/texts/2002/17$\_$1/93$\_$101.pdf

\bibitem[Ra3]{Ra3}
A. Ravsky, {\em The topological and algebraical properties of
paratopological groups}, Ph.D. Thesis. -- Lviv University, 2002
(in Ukrainian).

\bibitem[Ra4]{Ra4}
A. Ravsky, {\it Cone topologies of paratopological groups}.
http://arxiv.org/abs/1406.2993


\bibitem[St]{St}
M. Stone, {\em Applications of the theory of Boolean rings to general topology},
Trans. Amer. Math. Soc. {\bf 41} (1937), 375--481.

\bibitem[ST]{ST} M.~Sanchis, M.~Tkachenko,
{\em Totally Lindel\"of and totally $\omega$-narrow paratopological groups},
 Topology Appl. {\bf 155}:4 (2008), 322-334.

\bibitem[ST2]{ST2}
I.~S\'anchez, M. Tkachenko, {\em Products of bounded subsets of paratopological groups},
Topology Appl., (to appear).

\bibitem[XieLiTu]{XieLiTu}
Li-Hong Xie, Piyu Li, Jin-Ji Tu, {\em Notes on (regular)
$T_3$-reflections in the category of semitopological groups},
Topology Appl., {\bf 178} (2014), 46--55.
http://www.sciencedirect.com/science/article/pii/S0166864114003757

\end{thebibliography}
\end{document}